\documentclass[11pt]{amsart}
\usepackage{amsmath,amssymb,mathrsfs}
\usepackage{amsmath,amssymb}
\newtheorem{theorem}{Theorem}[section]
\newtheorem{proposition}[theorem]{Proposition}
\newtheorem{corollary}[theorem]{Corollary}
\newtheorem{lemma}[theorem]{Lemma}

\begin{document}

\title[Self-similar solutions to the Ricci flow]{Rotational symmetry of self-similar solutions to the Ricci flow}
\author{Simon Brendle}
\begin{abstract}
Let $(M,g)$ be a three-dimensional steady gradient Ricci soliton which is non-flat and $\kappa$-noncollapsed. We prove that $(M,g)$ is isometric to the Bryant soliton up to scaling. This solves a problem mentioned in Perelman's first paper \cite{Perelman1}.
\end{abstract}
\address{Department of Mathematics \\ Stanford University \\ Stanford, CA 94305}
\thanks{The author was supported in part by the National Science Foundation under grant DMS-0905628.}
\maketitle

\section{Introduction}

Self-similar solutions play a central role in the study of the Ricci flow, and have been studied extensively in connection with singularity formation; see e.g. the work of R.~Hamilton \cite{Hamilton-survey} and G.~Perelman \cite{Perelman1}, \cite{Perelman2}, \cite{Perelman3}. There are three basic types of self-similar solutions, which are referred to as shrinking solitons; steady solitons; and expanding solitons. A steady Ricci soliton $(M,g)$ is characterized by the fact that $2 \, \text{\rm Ric} = \mathscr{L}_X(g)$ for some vector field $X$. If the vector field $X$ is the gradient of a function, we say that $(M,g)$ is a steady gradient Ricci soliton.

The simplest example of a steady Ricci soliton is the cigar soliton in dimension $2$, which was found by Hamilton (cf. \cite{Hamilton-survey}). R.~Bryant \cite{Bryant} has discovered a steady Ricci soliton in dimension $3$, which is rotationally symmetric. Moreover, Bryant showed that there are no other complete steady Ricci solitons in dimension $3$ which are rotationally symmetric. While additional examples are known in higher dimensions (see e.g. \cite{Ivey}), the Bryant soliton is so far the only known example of a non-flat steady Ricci soliton in dimension $3$. It is an interesting question whether any three-dimensional steady Ricci soliton is necessarily rotationally symmetric. Perelman mentions the uniqueness problem for steady Ricci solitons in his first paper (see \cite{Perelman1}, page 32, lines 8-9), without however indicating a strategy for a possible proof. 

In this paper, we prove the uniqueness of the Bryant soliton under a noncollapsing assumption, as proposed by Perelman:

\begin{theorem}
\label{main.theorem}
Let $(M,g)$ be a three-dimensional complete steady gradient Ricci soliton which is non-flat and $\kappa$-noncollapsed. Then $(M,g)$ is rotationally symmetric, and is therefore isometric to the Bryant soliton up to scaling. 
\end{theorem}

We note that several authors have obtained uniqueness results for the Bryant soliton and its higher dimensional counterparts under various additional assumptions. We refer to \cite{Cao-survey}, \cite{Cao-Chen}, \cite{Cao-Catino-Chen-Mantegazza-Mazzieri}, and \cite{Chen-Wang} for details.

We now outline the main steps involved in the proof of Theorem \ref{main.theorem}. Let $(M,g)$ be a three-dimensional complete steady gradient Ricci soliton which is non-flat and $\kappa$-noncollapsed. We may write $\text{\rm Ric} = D^2 f$ for some real-valued function $f$. For abbreviation, we put $X = \nabla f$. Moreover, we denote by $\Phi_t$ the one-parameter group of diffeomorphisms generated by the vector field $-X$. We may assume without loss of generality that $R + |\nabla f|^2 = 1$.

In Section \ref{asymptotics}, we analyze the asymptotic geometry of $(M,g)$. The local version of the Hamilton-Ivey pinching estimate established by B.L.~Chen \cite{Chen} implies that $(M,g)$ has positive sectional curvature. It then follows from work of Perelman \cite{Perelman1} that the flow $(M,g(t))$ is asymptotic to a family of shrinking cylinders near infinity. This fact plays a fundamental role in our analysis. We next show that the restriction of the scalar curvature to the level surface $\{f=r\}$ satisfies $R = \frac{1}{r} + O(r^{-\frac{5}{4}})$. As a consequence, the intrinsic Gaussian curvature of the level surface $\{f=r\}$ equals $\frac{1}{2r} + O(r^{-\frac{5}{4}})$. This can be viewed as a refined roundness estimate for the level surface $\{f=r\}$. 

In Section \ref{approx.killing.vector.fields}, we construct a collection of approximate Killing vector fields near infinity. More precisely, we construct three vector fields $U_1,U_2,U_3$ such that $|\mathscr{L}_{U_a}(g)| \leq O(r^{-\frac{1}{8}})$ and $|\Delta U_a + D_X U_a| \leq O(r^{-\frac{9}{16}})$. Moreover, we show that the vector fields $U_1,U_2,U_3$ satisfy 
\[\sum_{a=1}^3 U_a \otimes U_a = r \, (e_1 \otimes e_1 + e_2 \otimes e_2 + O(r^{-\frac{1}{8}})),\] 
where $\{e_1,e_2\}$ is a local orthonormal frame on the level set $\{f=r\}$. 

In Section \ref{calculation}, we consider a vector field $W$ which satisfies the elliptic equation $\Delta W + D_X W = 0$. We then consider the Lie derivative $h = \mathscr{L}_W(g)$. This tensor turns out to satisfy the equation 
\begin{equation} 
\label{lich.eq}
\Delta_L h + \mathscr{L}_X(h) = 0. 
\end{equation}
Here, $\Delta_L$ denotes the Lichnerowicz Laplacian; that is, 
\[\Delta_L h_{ik} = \Delta h_{ik} + 2 \, R_{ijkl} \, h^{jl} - \text{\rm Ric}_i^l \, h_{kl} - \text{\rm Ric}_k^l \, h_{il}.\] 

In Section \ref{pde.for.vector.fields}, we assume that a vector field $Q$ satisfying $|Q| \leq O(r^{-\frac{1}{2}-2\varepsilon})$ is given. We then construct a vector field $V$ such that $\Delta V + D_X V = Q$ and $|V| \leq O(r^{-\frac{1}{2}-\varepsilon})$. In order to construct the vector field $V$, we solve the Dirichlet problem on a sequence of domains which exhaust $M$. In order to be able to pass to the limit, we need uniform estimates for solutions of the equation $\Delta V + D_X V = Q$. These estimates are established using a delicate blow-down analysis; see Proposition \ref{uniform.bounds} below. 

In Section \ref{analysis.of.lich.eq}, we consider a symmetric $(0,2)$-tensor $h$ which solves the equation (\ref{lich.eq}) and satisfies $|h| \leq O(r^{-\varepsilon})$ at infinity. Note that such a tensor $h$ need not vanish identically. Indeed, the Ricci tensor of $(M,g)$ is a non-trivial solution of the equation (\ref{lich.eq}), which falls off like $r^{-1}$ at infinity. However, we are able to show that any solution of (\ref{lich.eq}) with $|h| \leq O(r^{-\varepsilon})$ is of the form $h = \lambda \, \text{\rm Ric}$ for some constant $\lambda \in \mathbb{R}$; see Theorem \ref{lichnerowicz.equation} below. The proof of Theorem \ref{lichnerowicz.equation} again relies on a parabolic blow-down argument. We also use an inequality due to G.~Anderson and B.~Chow \cite{Anderson-Chow} for solutions of the parabolic Lichnerowicz equation. Related ideas were used in earlier work of M.~Gursky \cite{Gursky} and R.~Hamilton \cite{Hamilton}.

Finally, in Section \ref{sym.prin}, we establish a crucial symmetry principle. To explain this, suppose that $U$ is a vector field on $(M,g)$ such that $|\mathscr{L}_U(g)| \leq O(r^{-2\varepsilon})$ and $|\Delta U + D_X U| \leq O(r^{-\frac{1}{2}-2\varepsilon})$ for some small constant $\varepsilon>0$. Using the results in Section \ref{pde.for.vector.fields}, we can find a vector field $V$ such that $\Delta V + D_X V = \Delta U + D_X U$ and $|V| \leq O(r^{-\frac{1}{2}-\varepsilon})$. Therefore, the vector field $W = U - V$ satisfies $\Delta W + D_X W = 0$. Consequently, the Lie derivative $h = \mathscr{L}_W(g)$ is a solution of the equation (\ref{lich.eq}). Moreover, we show that $|h| \leq O(r^{-\varepsilon})$ at infinity. Thus, $h = \lambda \, \text{\rm Ric}$ for some constant $\lambda \in \mathbb{R}$. From this, we deduce that the vector field $\hat{U} := W - \frac{1}{2} \, \lambda \, X$ is a Killing vector field. Moreover, the Killing vector field $\hat{U}$ agrees with the original vector field $U$ up to terms of order $O(r^{\frac{1}{2}-\varepsilon})$. 

Applying this symmetry principle to the approximate Killing vector fields $U_1,U_2,U_3$ constructed in Section \ref{approx.killing.vector.fields}, we obtain three exact Killing vector fields $\hat{U}_1,\hat{U}_2,\hat{U}_3$ on $(M,g)$ with the property that $\langle \hat{U}_a,X \rangle = 0$ and 
\[\sum_{a=1}^3 \hat{U}_a \otimes \hat{U}_a = r \, (e_1 \otimes e_1 + e_2 \otimes e_2 + O(r^{-\varepsilon})),\] 
where $\{e_1,e_2\}$ is a local orthonormal frame on the level surface $\{f=r\}$. In particular, at each point sufficiently far out at infinity, the span of the vector fields $\hat{U}_1,\hat{U}_2,\hat{U}_3$ is two-dimensional.

Finally, let us mention some related results. Our method of proof is inspired in part by the beautiful work of L.~Simon and B.~Solomon on the uniqueness of minimal hypersurfaces in $\mathbb{R}^{n+1}$ which are asymptotic to a given cone at infinity (cf. \cite{Simon}, \cite{Simon-Solomon}). X.J.~Wang \cite{Wang} has obtained a uniqueness theorem for convex translating solutions to the mean curvature flow in $\mathbb{R}^3$. The argument in \cite{Wang} is quite different from ours and relies in a crucial way on a classical theorem of Bernstein (cf. \cite{Hopf}). Finally, the uniqueness problem for the Bryant soliton shares some common features with the black hole uniqueness theorems in general relativity (see e.g. \cite{Hawking-Ellis}, \cite{Ionescu-Klainerman}). 

It is a pleasure to thank Professors Huai-Dong Cao, Gerhard Huisken, Sergiu Klainerman, Leon Simon, Brian White, for discussions. The author is grateful to Meng Zhu for comments on an earlier version of this paper.

\section{The asymptotic geometry of $(M,g)$}

\label{asymptotics}

Throughout this paper, we assume that $(M,g)$ is a three-dimensional complete steady gradient Ricci soliton which is $\kappa$-noncollapsed and non-flat. It follows from Theorem 1.3 in \cite{Zhang} that $(M,g)$ has positive scalar curvature (see also \cite{Cao-Chen}, Proposition 2.2). It is well known that the sum $R + |\nabla f|^2$ is constant. By scaling, we may assume that $R + |\nabla f|^2 = 1$. Since $R \geq 0$, it follows that $|\nabla f|^2 \leq 1$. Hence, if we denote by $\Phi_t$ the flow generated by the vector field $-X$, then $\Phi_t$ is defined for all $t \in \mathbb{R}$, and the metrics $\Phi_t^*(g)$ evolve by the Ricci flow. 

\begin{proposition}
\label{positive.sectional.curvature}
The manifold $(M,g)$ has bounded curvature, and the sectional curvature is strictly positive.
\end{proposition} 

\textbf{Proof.} 
It follows from a result of Chen that $(M,g)$ has nonnegative sectional curvature (see \cite{Chen}, Corollary 2.4). Since $R + |\nabla f|^2 \leq 1$, we conclude that $(M,g)$ has bounded curvature. It remains to show that $(M,g)$ has positive sectional curvature. Suppose this is false. Then the manifold $(M,g)$ locally splits as a product, and the universal cover of $(M,g)$ is isometric to the cigar soliton crossed with a line. This contradicts our assumption that $(M,g)$ is $\kappa$-noncollapsed. \\

We next analyze the asymptotic geometry of $(M,g)$ near infinity. We will frequently use the identity 
\begin{equation} 
\label{evol.of.scalar.curvature}
-\langle X,\nabla R \rangle = \Delta R + 2 \, |\text{\rm Ric}|^2. 
\end{equation} 
This identity is a consequence of the evolution equation for the scalar curvature under the Ricci flow (cf. \cite{Brendle-book}, Section 2.4).

The following result is a direct consequence of Perelman's compactness theorem for ancient $\kappa$-solutions:

\begin{proposition}[G.~Perelman \cite{Perelman1}, \cite{Perelman2}]
\label{basics}
Let $p_m$ be a sequence of points going to infinity. Then $|\langle X,\nabla R \rangle| \leq O(1) \, R^2$ at the point $p_m$. Moreover, if $d(p_0,p_m)^2 \, R(p_m) \to \infty$, then we have $|\nabla R| \leq o(1) \, R^{\frac{3}{2}}$ and $|\langle X,\nabla R \rangle + R^2| \leq o(1) \, R^2$ at the point $p_m$.
\end{proposition}

\textbf{Proof.} 
It follows from results in Section 1.5 of \cite{Perelman2} that $|\Delta R| \leq O(1) \, R^2$. Using (\ref{evol.of.scalar.curvature}), we conclude that $|\langle X,\nabla R \rangle| \leq O(1) \, R^2$. This proves the first statement. 

We now describe the proof of the second statement. To that end, we assume that $d(p_0,p_m) \, R(p_m)^2 \to \infty$. Let us consider the rescaled flows 
\[\hat{g}^{(m)}(t) = r_m^{-1} \, \Phi_{r_m t}^*(g),\] 
where $r_m = R(p_m)^{-1}$. It follows from Perelman's compactness theorem for ancient $\kappa$-solutions that the flows $(M,\hat{g}^{(m)}(t),p_m)$, $t \in (-\infty,0]$, converge in the Cheeger-Gromov sense to a non-flat ancient $\kappa$-solution $(\overline{M},\overline{g}(t))$, $t \in (-\infty,0]$ (see \cite{Perelman1}, Theorem 11.7). By Theorem 5.35 in \cite{Morgan-Tian}, the manifold $(\overline{M},\overline{g}(0))$ splits off a line. By the strict maximum principle, the limit flow $(\overline{M},\overline{g}(t))$, $t \in (-\infty,0]$, is isometric to a product of a two-dimensional ancient $\kappa$-solution with a line. By Theorem 11.3 in \cite{Perelman1}, the universal cover of $(\overline{M},\overline{g}(t))$ is a round cylinder for each $t \in (-\infty,0]$. From this, we deduce that $|\nabla R| \leq o(1) \, R^{\frac{3}{2}}$, $|\Delta R| \leq o(1) \, R^2$, and $2 \, |\text{\rm Ric}|^2 = (1 + o(1)) \, R^2$ at the point $p_m$. Using (\ref{evol.of.scalar.curvature}), we conclude that $-\langle X,\nabla R \rangle = \Delta R + 2 \, |\text{\rm Ric}|^2 = (1+o(1)) \, R^2$. \\

\begin{corollary}
\label{scalar.curvature.goes.to.zero}
The scalar curvature converges to $0$ at infinity.
\end{corollary} 

\textbf{Proof.} 
Suppose this is false. Then we can find a sequence of points $p_m$ going to infinity such that $\liminf_{m \to \infty} R(p_m) > 0$. Using Proposition \ref{basics}, we obtain $|\langle X,\nabla R \rangle + R^2| \leq o(1)$ and $|\nabla R| \leq o(1)$ at the point $p_m$. Since $|X| \leq 1$, it follows that $|\langle X,\nabla R \rangle| \leq o(1)$ at the point $p_m$. Putting these facts together, we conclude that $R(p_m) = o(1)$, contrary to our assumption. \\

By Corollary \ref{scalar.curvature.goes.to.zero}, we can find a point $p_0 \in M$ such that $R(p_0) = \sup_M R$. At the point $p_0$, we have 
\[0 = \partial_i R = -2 \, D_{i,j}^2 f \, \partial^j f.\] By Proposition \ref{positive.sectional.curvature}, the Hessian of $f$ is positive definite at each point in $M$. Consequently, the point $p_0$ is a critical point of $f$. Moreover, we can find positive constants $c_1$ and $c_2$ such that 
\[c_1 \, d(p_0,p) \leq f(p) \leq c_2 \, d(p_0,p)\] 
outside of a compact set (see also \cite{Cao-Chen}, Proposition 2.3). Without loss of generality, we may assume that $\inf_M f \geq 1$.

\begin{proposition}[H.~Guo \cite{Guo}]
\label{guo}
The scalar curvature satisfies $f \, R = 1+o(1)$ as $p \to \infty$.
\end{proposition}

\textbf{Proof.} 
Using Corollary \ref{scalar.curvature.goes.to.zero} and the identity $R + |\nabla f|^2 = 1$, we obtain $|\nabla f|^2 \to 1$ as $p \to \infty$. In particular, we have $|\nabla f|^2 \geq \frac{1}{2}$ outside a compact set. Using Proposition \ref{basics}, we obtain $-\langle X,\nabla R \rangle \leq C \, R^2$, hence 
\[\Big \langle X,\nabla \Big ( \frac{1}{R} - 2C \, f \Big ) \Big \rangle \leq C \, (1 - 2 \, |\nabla f|^2) \leq 0\] 
outside a compact set. Integrating this inequality along the integral curves of $X$ gives 
\[\sup_M \Big ( \frac{1}{R} - 2C \, f \Big ) < \infty.\] 
Consequently, $\inf_M f \, R > 0$. In particular, we have $d(p_0,p)^2 \, R(p) \to \infty$ at infinity. Using Proposition \ref{basics} again, we conclude that 
\[|\langle X,\nabla R \rangle + R^2| \leq o(1) \, R^2\] 
near infinity. Since $1 - |\nabla f|^2 = R \to 0$ at infinity, we conclude that 
\[\Big \langle X,\nabla \Big ( \frac{1}{R} - f \Big ) \Big \rangle = 1 - |\nabla f|^2 - \frac{1}{R^2} \, (\langle X,\nabla R \rangle + R^2) = o(1).\] 
Integrating this inequality along the integral curves of $X$, we obtain 
\[\frac{1}{R} = (1+o(1)) \, f,\] 
as claimed. \\

Using work of Perelman \cite{Perelman1}, we can determine the asymptotic geometry of $(M,g)$ near infinity:

\begin{proposition}[cf. \cite{Perelman1}]
Let $p_m$ be a sequence of marked points going to infinity. Consider the rescaled metrics 
\[\hat{g}^{(m)}(t) = r_m^{-1} \, \Phi_{r_m t}^*(g),\] 
where $r_m = f(p_m)$. As $m \to \infty$, the flows $(M,\hat{g}^{(m)}(t),p_m)$ converge in the Cheeger-Gromov sense to a family of shrinking cylinders $(S^2 \times \mathbb{R},\overline{g}(t))$, $t \in (0,1)$. The metric $\overline{g}(t)$ is given by 
\begin{equation} 
\label{shrinking.cylinders}
\overline{g}(t) = (2-2t) \, g_{S^2} + dz \otimes dz, 
\end{equation}
where $g_{S^2}$ denotes the standard metric on $S^2$ with constant Gaussian curvature $1$. Furthermore, the rescaled vector fields $r_m^{\frac{1}{2}} \, X$ converge in $C_{loc}^\infty$ to the axial vector field $\frac{\partial}{\partial z}$ on $S^2 \times \mathbb{R}$.
\end{proposition}

\textbf{Proof.} 
It follows from Proposition \ref{guo} that the flows $(M,\hat{g}^{(m)}(t),p_m)$, $t \in (-\infty,1)$, converge in the Cheeger-Gromov sense to a non-flat ancient $\kappa$-solution $(\overline{M},\overline{g}(t))$, $t \in (-\infty,1)$. By Theorem 5.35 in \cite{Morgan-Tian}, the limit flow $(\overline{M},\overline{g}(t))$ is isometric to a product of a two-dimensional ancient $\kappa$-solution with a line (see \cite{Perelman1}, Theorem 11.7). Note that $M$ is homeomorphic to $\mathbb{R}^3$ and in particular does not contain an embedded $\mathbb{RP}^2$. Consequently, $\overline{M}$ cannot contain an embedded $\mathbb{RP}^2$. By Theorem 11.3 in \cite{Perelman1}, we conclude that $(\overline{M},\overline{g}(t))$ is a family of round cylinders, i.e. $\overline{M} = S^2 \times \mathbb{R}$ and $\overline{g}(t) = (2-2t) \, g_{S^2} + dz \otimes dz$ for each $t \in (-\infty,1)$. 

It remains to analyze the limit of the rescaled vector fields $\hat{X}^{(m)} = r_m^{\frac{1}{2}} \, X$. Using the identity $1 - |X| = O(r^{-1})$, we obtain 
\[\limsup_{m \to \infty} \sup_{\{r_m - \delta^{-1} \, \sqrt{r_m} \leq f \leq r_m + \delta^{-1} \, \sqrt{r_m}\}} \big | 1 - |\hat{X}^{(m)}|_{\hat{g}^{(m)}(0)} \big | = 0\] 
for any given $\delta \in (0,1)$. Moreover, we have $|D^l X| \leq C \, |D^{l-1} \text{\rm Ric}| = O(r^{-\frac{l+1}{2}})$ for all $l \geq 1$. This implies 
\[\limsup_{m \to \infty} \sup_{\{r_m - \delta^{-1} \, \sqrt{r_m} \leq f \leq r_m + \delta^{-1} \, \sqrt{r_m}\}} |D_{\hat{g}^{(m)}(0)}^l \hat{X}^{(m)}|_{\hat{g}^{(m)}(0)} = 0\] 
for any given $\delta \in (0,1)$ and $l \geq 1$. Hence, after passing to a subsequence, the vector fields $\hat{X}^{(m)}$ converge in $C_{loc}^\infty$ to a vector field $\overline{X}$ on the limit manifold $(S^2 \times \mathbb{R},\overline{g}(0))$. The limiting vector field $\overline{X}$ is parallel with respect to the metric $\overline{g}(0)$, and we have $|\overline{X}|_{\overline{g}(0)} = 1$. Thus, $\overline{X}$ can be identified with the axial vector field $\frac{\partial}{\partial z}$ on $S^2 \times \mathbb{R}$. \\

In the remainder of this section, we establish a roundness estimate for the level surfaces $\{f=r\}$. The proof of this estimate requires several lemmata.

\begin{lemma} 
\label{grad.of.scalar.curvature}
On the level surface $\{f=r\}$, we have 
\[2 \, \text{\rm Ric}(\nabla f,\nabla f) = -\langle X,\nabla R \rangle = O(r^{-2}).\] 
\end{lemma} 

\textbf{Proof.} 
The identity (\ref{evol.of.scalar.curvature}) implies that $2 \, \text{\rm Ric}(\nabla f,\nabla f) = -\langle X,\nabla R \rangle = \Delta R + 2 \, |\text{\rm Ric}|^2 = O(r^{-2})$. \\

\begin{lemma} 
\label{mean.curvature}
The mean curvature of the level surface $\{f=r\}$ equals $\frac{1 + o(1)}{r}$.
\end{lemma} 

\textbf{Proof.} 
The mean curvature of the level surface $\{f=r\}$ is given by 
\[H = \frac{1}{|\nabla f|} \, R - \frac{1}{|\nabla f|^3} \, \text{\rm Ric}(\nabla f,\nabla f).\] 
Hence, the assertion follows from Proposition \ref{guo} and Lemma \ref{grad.of.scalar.curvature}. \\

\begin{lemma}
\label{ricci.tensor}
The tensor $T = 2 \, \text{\rm Ric} - R \, g + R \, df \otimes df$ satisfies $|T| \leq O(r^{-\frac{3}{2}})$ and $|DT| \leq O(r^{-2})$.
\end{lemma} 

\textbf{Proof.} 
In dimension $3$, the Riemann curvature tensor can be written in the form 
\begin{align*} 
R_{ijkl} 
&= \text{\rm Ric}_{ik} \, g_{jl} - \text{\rm Ric}_{il} \, g_{jk} - \text{\rm Ric}_{jk} \, g_{il} + \text{\rm Ric}_{jl} \, g_{ik} \\ 
&- \frac{1}{2} \, R \, (g_{ik} \, g_{jl} - g_{il} \, g_{jk}). 
\end{align*}
This implies 
\begin{align*} 
D_i \text{\rm Ric}_{jk} - D_j \text{\rm Ric}_{ik} 
&= R_{ijkl} \, D^l f \\ 
&= \text{\rm Ric}_{ik} \, D_j f - \text{\rm Ric}_{jk} \, D_i f \\ 
&- \frac{1}{2} \, (D_j R + R \, D_j f) \, g_{ik} + \frac{1}{2} \, (D_i R + R \, D_i f) \, g_{jk}, 
\end{align*} 
hence 
\begin{align} 
\label{T}
&2 \, (D_i \text{\rm Ric}_{jk} - D_j \text{\rm Ric}_{ik}) \, D^j f \notag \\ 
&= T_{ik} \, |\nabla f|^2 - \langle \nabla R,\nabla f \rangle \, g_{ik} + R^2 \, D_i f \, D_k f \\ 
&+ D_i R \, D_k f + D_k R \, D_i f. \notag 
\end{align} 
By Shi's estimate, the covariant derivatives of the curvature tensor are bounded by $O(r^{-\frac{3}{2}})$. Consequently, the identity (\ref{T}) implies that $|T| \leq O(r^{-\frac{3}{2}})$. Moreover, if we differentiate (\ref{T}), we obtain $|DT| \leq O(r^{-2})$. \\

\begin{lemma} 
\label{evolution}
We have 
\[|\langle X,\nabla R \rangle + \Delta_\Sigma R + R^2| \leq O(r^{-\frac{5}{2}}),\] 
where $\Delta_\Sigma$ denotes the Laplacian on the level surface $\{f=r\}$.
\end{lemma} 

\textbf{Proof.} 
Differentiating the identity (\ref{evol.of.scalar.curvature}), we obtain 
\[-(D^2 R)(X,X) - \langle D_X X,\nabla R \rangle = \langle X,\nabla (\Delta R + 2 \, |\text{\rm Ric}|^2) \rangle.\] 
Since $\nabla R = -2 \, D_X X$, it follows that 
\[-(D^2 R)(X,X) =  -\frac{1}{2} \, |\nabla R|^2 + \langle X,\nabla (\Delta R + 2 \, |\text{\rm Ric}|^2) \rangle.\] 
Using Shi's estimates, we obtain $|\nabla R|^2 \leq O(r^{-3})$ and $|\nabla (\Delta R + 2 \, |\text{\rm Ric}|^2)| \leq O(r^{-\frac{5}{2}})$. Consequently, we have 
\begin{equation} 
\label{radial.second.derivative}
|(D^2 R)(X,X)| \leq O(r^{-\frac{5}{2}}). 
\end{equation}
Moreover, it follows from Lemma \ref{grad.of.scalar.curvature} and Lemma \ref{mean.curvature} that  
\begin{equation} 
\label{mean.curvature.term}
|H \, \langle X,\nabla R \rangle| \leq O(r^{-3}). 
\end{equation} 
Combining (\ref{radial.second.derivative}) and (\ref{mean.curvature.term}) gives 
\[|\Delta R - \Delta_\Sigma R| \leq O(r^{-\frac{5}{2}}).\] 
Combining this inequality with (\ref{evol.of.scalar.curvature}), we obtain 
\[\big | \Delta_\Sigma R + \langle X,\nabla R \rangle + 2 \, |\text{\rm Ric}|^2 \big | \leq O(r^{-\frac{5}{2}}).\] 
On the other hand, it follows from Lemma \ref{ricci.tensor} that 
\[2 \, |\text{\rm Ric}| = |R \, (g - df \otimes df)| + O(r^{-\frac{3}{2}}) = \sqrt{2} \, R + O(r^{-\frac{3}{2}}).\] 
Putting these facts together, the assertion follows. \\

We next establish a Poincar\'e-type inequality for the restriction of the scalar curvature to a level surface $\{f=r\}$. Our argument uses the Kazdan-Warner identity (cf. \cite{Kazdan-Warner}), and is inspired in part by work of M.~Struwe on the Calabi flow on the two-sphere (cf. \cite{Struwe}, p.~263). In the sequel, we denote by $\mu(r)$ the mean value of the scalar curvature over the level surface $\{f=r\}$, so that 
\[\int_{\{f=r\}} (R - \mu(r)) = 0.\] 
Note that $\mu(r) = \frac{1+o(1)}{r}$ by Proposition \ref{guo}. 

\begin{lemma}
\label{eigenvalue.estimate}
We have 
\[\int_{\{f=r\}} |\nabla^\Sigma R|^2 \geq \frac{2}{r} \, \bigg ( \int_{\{f=r\}} (R-\mu(r))^2 \bigg ) - O(r^{-4})\] 
if $r$ is sufficiently large.
\end{lemma}

\textbf{Proof.} 
Let us fix $r$ sufficiently large. Let $0 = \nu_0 < \nu_1 \leq \nu_2 \leq \nu_3 \leq \hdots$ denote the eigenvalues of the Laplace operator on the level surface $\{f=r\}$, and let $\psi_0,\psi_1,\psi_2,\psi_3,\hdots$ denote the associated eigenfunctions. We assume that the eigenfunctions are normalized so that $\int_{\{f=r\}} \psi_j^2 = 1$ for each $j$. When $r$ is large, the surface $\{f=r\}$ equipped with the rescaled metric $\frac{1}{2r} \, g$ is $C^\infty$ close to the standard two-sphere with constant Gaussian curvature $1$. Consequently, $\nu_1 = \frac{1+o(1)}{r}$, $\nu_2 = \frac{1+o(1)}{r}$, $\nu_3 = \frac{1+o(1)}{r}$, and $\nu_4 = \frac{3+o(1)}{r}$. 

Let $K$ denote the intrinsic Gaussian curvature of the level surface $\{f=r\}$. Using the Gauss equations, we obtain 
\[R - \frac{2}{|\nabla f|^2} \, \text{\rm Ric}(\nabla f,\nabla f) = 2 \, R(e_1,e_2,e_1,e_2) = 2K + O(r^{-2}).\] 
Using Lemma \ref{grad.of.scalar.curvature}, we conclude that $|2K - R| \leq O(r^{-2})$, hence 
\[\bigg ( \int_{\{f=r\}} (2K - R)^2 \bigg )^{\frac{1}{2}} \leq O(r^{-\frac{3}{2}}).\] 
On the other hand, it follows from the Kazdan-Warner identity (see \cite{Kazdan-Warner}, Theorem 8.8) that 
\[\sum_{j=1}^3 \bigg | \int_{\{f=r\}} (2K-\mu(r)) \, \psi_j \bigg | \leq o(1) \, \bigg ( \int_{\{f=r\}} (2K-\mu(r))^2 \bigg )^{\frac{1}{2}}\] 
Putting these facts together, we obtain 
\[\sum_{j=1}^3 \bigg | \int_{\{f=r\}} (R-\mu(r)) \, \psi_j \bigg | \leq o(1) \, \bigg ( \int_{\{f=r\}} (R-\mu(r))^2 \bigg )^{\frac{1}{2}} + O(r^{-\frac{3}{2}}).\] 
Thus, we conclude that 
\begin{align*} 
&\int_{\{f=r\}} |\nabla^\Sigma R|^2 - \nu_4 \int_{\{f=r\}} (R-\mu(r))^2 \\ 
&= \sum_{j=1}^\infty (\nu_j - \nu_4) \, \bigg ( \int_{\{f=r\}} (R-\mu(r)) \, \psi_j \bigg )^2 \\ 
&\geq -\frac{3}{r} \, \sum_{j=1}^3 \bigg ( \int_{\{f=r\}} (R-\mu(r)) \, \psi_j \bigg )^2 \\ 
&\geq -o(r^{-1}) \, \bigg ( \int_{\{f=r\}} (R-\mu(r))^2 \bigg ) - O(r^{-4}). 
\end{align*} 
Since $\nu_4 = \frac{3+o(1)}{r}$, the assertion follows. \\

We now prove an important roundness estimate. 

\begin{proposition} 
\label{roundness.1}
We have 
\[\int_{\{f=r\}} (R - \mu(r))^2 \leq O(r^{-2})\] 
if $r$ is sufficiently large.
\end{proposition} 

\textbf{Proof.} 
By definition of $\mu(r)$, we have $\int_{\{f=r\}} (R-\mu(r)) = 0$. This implies 
\begin{align*} 
&\frac{d}{dr} \bigg ( \int_{\{f=r\}} (R - \mu(r))^2 \bigg ) \\ 
&= 2 \int_{\{f=r\}} (R - \mu(r)) \, \Big ( \frac{\langle X,\nabla R \rangle}{|X|^2} - \mu'(r) \Big ) + \int_{\{f=r\}} \frac{H}{|X|} \, (R - \mu(r))^2 \\ 
&= 2 \int_{\{f=r\}} (R - \mu(r)) \, \Big ( \frac{\langle X,\nabla R \rangle}{|X|^2} + \mu(r)^2 \Big ) + \int_{\{f=r\}} \frac{H}{|X|} \, (R - \mu(r))^2 \\ 
&= 2 \int_{\{f=r\}} |\nabla^\Sigma R|^2 - \int_{\{f=r\}} \Big ( 2R + 2\mu(r) - \frac{H}{|X|} \Big ) \, (R - \mu(r))^2 \\ 
&+ 2 \int_{\{f=r\}} (R - \mu(r)) \, \Big ( \frac{\langle X,\nabla R \rangle}{|X|^2} + \Delta_\Sigma R + R^2 \Big ). 
\end{align*} 
It follows from Lemma \ref{eigenvalue.estimate} that 
\[\int_{\{f=r\}} |\nabla^\Sigma R|^2 \geq \frac{2}{r} \, \bigg ( \int_{\{f=r\}} (R-\mu(r))^2 \bigg ) - O(r^{-4}).\] 
Moreover, we have $2R+2\mu(r)-\frac{H}{|X|} = \frac{3+o(1)}{r}$. Finally, we have 
\[\bigg | \frac{\langle X,\nabla R \rangle}{|X|^2} + \Delta_\Sigma R + R^2 \bigg | \leq O(r^{-\frac{5}{2}})\] 
by Lemma \ref{evolution}. Putting these facts together, we obtain 
\begin{align*} 
\frac{d}{dr} \bigg ( \int_{\{f=r\}} (R - \mu(r))^2 \bigg ) 
&\geq \frac{1-o(1)}{r} \int_{\{f=r\}} (R - \mu(r))^2 \\ 
&- O(r^{-\frac{5}{2}}) \, \int_{\{f=r\}} |R - \mu(r)| \\ 
&- O(r^{-4}). 
\end{align*} 
Using Young's inequality, we conclude that 
\begin{align*} 
\frac{d}{dr} \bigg ( \int_{\{f=r\}} (R - \mu(r))^2 \bigg ) 
&\geq -O(r^{-4}) \, \text{\rm vol}(\{f=r\}) - O(r^{-4}) \\ 
&\geq -O(r^{-3}). 
\end{align*}
Clearly, 
\[\int_{\{f=r\}} (R - \mu(r))^2 \to 0\] 
as $r \to \infty$. Putting these facts together, we obtain 
\[\int_{\{f=r\}} (R-\mu(r))^2 \leq O(r^{-2}),\] 
as claimed. \\

\begin{corollary} 
\label{roundness.2}
We have 
\begin{align*} 
&\sup_{\{f=r\}} |R - \mu(r)| \leq O(r^{-\frac{5}{4}}), \\ 
&\sup_{\{f=r\}} |\nabla^\Sigma R| \leq O(r^{-\frac{7}{4}}), \\ 
&\sup_{\{f=r\}} |\Delta_\Sigma R| \leq O(r^{-\frac{9}{4}}).
\end{align*}
\end{corollary}

\textbf{Proof.} 
By Proposition \ref{roundness.1}, we have 
\[\int_{\{f=r\}} (R - \mu(r))^2 \leq O(r^{-2}).\] 
Moreover, it follows from Shi's estimates that 
\[\sup_{\{f=r\}} |D_\Sigma^l R| \leq O(r^{-\frac{l+2}{2}}).\] 
Hence, the assertion follows from standard interpolation inequalities (see e.g. \cite{Hamilton}, Corollary 12.7). \\

With the aid of Corollary \ref{roundness.2}, we can improve Proposition \ref{guo} as follows: 

\begin{proposition} 
\label{scalar.2}
We have $|\nabla R| \leq O(r^{-\frac{7}{4}})$ and $f \, R = 1 + O(r^{-\frac{1}{4}})$.
\end{proposition}

\textbf{Proof.} 
Using the estimates $|\nabla^\Sigma R| \leq O(r^{-\frac{7}{4}})$ and $|\langle X,\nabla R \rangle| \leq O(r^{-2})$, we obtain $|\nabla R| \leq O(r^{-\frac{7}{4}})$. This proves the first statement.

We now describe the proof of the second statement. By Corollary \ref{roundness.2}, we have $|\Delta_\Sigma R| \leq O(r^{-\frac{9}{4}})$. Hence, Lemma \ref{evolution} implies  
\[|\langle X,\nabla R \rangle + R^2| \leq O(r^{-\frac{9}{4}}).\] From this, we deduce that 
\[\Big \langle X,\nabla \Big ( \frac{1}{R} - f \Big ) \Big \rangle = 1 - |\nabla f|^2 - \frac{1}{R^2} \, (\langle X,\nabla R \rangle + R^2) = O(r^{-\frac{1}{4}}).\] 
Integrating this relation along the integral curves of $X$ gives 
\[\frac{1}{R} - f = O(r^{\frac{3}{4}}).\] From this, the assertion follows. \\

\begin{corollary} 
\label{geometry.of.level.sets}
The principal curvatures of the level surface $\{f=r\}$ are given by $\frac{1}{2r} + O(r^{-\frac{5}{4}})$. Moreover, the intrinsic Gaussian curvature of the level surface $\{f=r\}$ is given by $\frac{1}{2r} + O(r^{-\frac{5}{4}})$.
\end{corollary}

\section{Existence of approximate Killing vector fields near infinity}

\label{approx.killing.vector.fields}

In this section, we shall construct a collection of approximate Killing vector fields near infinity. It is easy to see that the level surfaces of $f$ are diffeomorphic to $S^2$. Hence, we can find a family of diffeomorphisms $F_r: S^2 \to \{f=r\} \subset M$ such that $\frac{\partial}{\partial r} F_r = \frac{X}{|X|^2}$. We define a metric $\gamma_r$ on $S^2$ by $\gamma_r = \frac{1}{2r} \, F_r^*(g)$. 

\begin{proposition} 
\label{induced.metrics}
We have 
\[\Big \| \frac{d}{dr} \gamma_r \Big \|_{C^l(S^2,\gamma_r)} \leq O(r^{-\frac{9}{8}})\]
for each $l \geq 0$.
\end{proposition}

\textbf{Proof.} 
By Corollary \ref{geometry.of.level.sets}, the principal curvatures of the level surface $\{f=r\}$ are given by $\frac{1}{2r} + O(r^{-\frac{5}{4}})$. Moreover, the normal velocity of the flow $F_r: S^2 \to M$ is $\frac{1}{|X|} = 1 + O(r^{-1})$. This implies 
\[\sup_{S^2} \Big | \frac{d}{dr} F_r^*(g) - \frac{1}{r} \, F_r^*(g) \Big |_{F_r^*(g)} \leq O(r^{-\frac{5}{4}}).\] From this, we deduce that 
\begin{equation} 
\label{C0.bound}
\sup_{S^2} \Big | \frac{d}{dr} \gamma_r \Big |_{\gamma_r} \leq O(r^{-\frac{5}{4}}). 
\end{equation}
Using the estimate $\sup_{\{f=r\}} |D^l \text{\rm Ric}| \leq O(r^{-\frac{l+2}{2}})$, we conclude that the manifold $(S^2,\gamma_r)$ has bounded curvature, and all the derivatives of the curvature are bounded as well. Using the inequality 
\[\sup_{\{r-\sqrt{r} \leq f \leq r+\sqrt{r}\}} \Big | D^l \big ( \mathscr{L}_{\frac{X}{|X|^2}}(g) \big ) \Big | \leq O(r^{-\frac{l+2}{2}}),\] 
we obtain 
\[\Big \| F_r^* \big ( \mathscr{L}_{\frac{X}{|X|^2}}(g) \big ) \Big \|_{C^l(S^2,\gamma_r)} \leq O(1).\] 
Since 
\[\frac{d}{dr} \gamma_r + \frac{1}{r} \, \gamma_r = \frac{1}{2r} \, F_r^* \big ( \mathscr{L}_{\frac{X}{|X|^2}}(g) \big ),\] 
we conclude that 
\begin{equation} 
\label{Cl.bound}
\Big \| \frac{d}{dr} \gamma_r \Big \|_{C^l(S^2,\gamma_r)} \leq O(r^{-1}) 
\end{equation}
for each $l \geq 0$. Using (\ref{C0.bound}), (\ref{Cl.bound}), and standard interpolation inequalities, the assertion follows. \\

By Proposition A.5 in \cite{Brendle-book}, the metrics $\gamma_r$ converge in $C^\infty$ to a smooth metric $\overline{\gamma}$ as $r \to \infty$. By Corollary \ref{geometry.of.level.sets}, the Gaussian curvature of the metric $\gamma_r$ is $1 + O(r^{-\frac{1}{4}})$. Consequently, the limit metric $\overline{\gamma}$ must have constant Gaussian curvature $1$. Moreover, Proposition \ref{induced.metrics} implies that 
\begin{equation} 
\Big \| \frac{d}{dr} \gamma_r \Big \|_{C^l(S^2,\overline{\gamma})} \leq O(r^{-\frac{9}{8}}), 
\end{equation}
hence 
\begin{equation} 
\label{difference.of.metrics}
\|\gamma_r - \overline{\gamma}\|_{C^l(S^2,\overline{\gamma})} \leq O(r^{-\frac{1}{8}}) 
\end{equation}
for each $l \geq 0$. 

Let $\overline{U}_1,\overline{U}_2,\overline{U}_3$ be three Killing vector fields on the round sphere $(S^2,\overline{\gamma})$ such that  
\begin{equation} 
\label{normalization}
\sum_{a=1}^3 \overline{U}_a \otimes \overline{U}_a = \frac{1}{2} \, (\overline{e}_1 \otimes \overline{e}_1 + \overline{e}_2 \otimes \overline{e}_2), 
\end{equation}
where $\{\overline{e}_1,\overline{e}_2\}$ is a local orthonormal frame on $(S^2,\overline{\gamma})$. Using (\ref{difference.of.metrics}), we obtain 
\begin{equation} 
\label{approximate.killing}
\|\mathscr{L}_{\overline{U}_a}(\gamma_r)\|_{C^l(S^2,\gamma_r)} = \|\mathscr{L}_{\overline{U}_a}(\gamma_r - \overline{\gamma})\|_{C^l(S^2,\gamma_r)} \leq O(r^{-\frac{1}{8}}). 
\end{equation} 
We can find three vector fields $U_1,U_2,U_3$ on $M$ with the property that the vector field $U_a$ is tangential to the level set $\{f=r\}$, and $F_r^* U_a = \overline{U}_a$ for $r$ sufficiently large. Clearly, $[U_a,\frac{X}{|X|^2}] = 0$ outside a compact set. This implies 
\begin{equation} 
\label{lie.bracket}
[U_a,X] = U_a(|X|^2) \, \frac{X}{|X|^2} 
\end{equation}
outside a compact set. Since $U_a(|X|^2) = -\langle U_a,\nabla R \rangle = O(r^{-1})$, we conclude that $|[U_a,X]| \leq O(r^{-1})$. Moreover, the inequality $\|\overline{U}_a\|_{C^l(S^2,\gamma_r)} \leq O(1)$ gives $\sup_{\{f=r\}} |D_\Sigma^l U_a| \leq O(r^{-\frac{l-1}{2}})$ for each $l \geq 0$. Since $[U_a,\frac{X}{|X|^2}] = 0$, we conclude that $\sup_{\{f=r\}} |D^l U_a| \leq O(r^{-\frac{l-1}{2}})$ for each $l \geq 0$. \\

\begin{proposition}
\label{rotation.vector.fields}
The vector fields $U_1,U_2,U_3$ on $(M,g)$ satisfy $|\mathscr{L}_{U_a}(g)| \leq O(r^{-\frac{1}{8}})$ and $|\Delta U_a + D_X U_a| \leq O(r^{-\frac{9}{16}})$. Moreover, we have 
\[\sum_{a=1}^3 U_a \otimes U_a = r \, (e_1 \otimes e_1 + e_2 \otimes e_2 + O(r^{-\frac{1}{8}})),\] 
where $\{e_1,e_2\}$ is a local orthonormal frame on the level set $\{f=r\}$. 
\end{proposition}

\textbf{Proof.} 
Let $\{e_1,e_2\}$ be a local orthonormal frame on the level surface $\{f=r\}$. Using (\ref{approximate.killing}) and (\ref{lie.bracket}), we obtain 
\begin{align*} 
&\langle D_{e_i} U_a,e_j \rangle + \langle D_{e_j} U_a,e_i \rangle = O(r^{-\frac{1}{8}}), \\ 
&\langle D_X U_a,e_j \rangle + \langle D_{e_j} U_a,X \rangle = \langle D_{U_a} X,e_j \rangle - \langle U_a,D_{e_j} X \rangle - \langle [U_a,X],e_j \rangle = 0, \\ 
&\langle D_X U_a,X \rangle = \langle D_{U_a} X,X \rangle - \langle [U_a,X],X \rangle = -\frac{1}{2} \, U_a(|X|^2) = O(r^{-1}). 
\end{align*} 
Therefore, the tensor $h_a = \mathscr{L}_{U_a}(g)$ satisfies 
\[\sup_{\{f=r\}} |h_a| \leq O(r^{-\frac{1}{8}}).\] 
Moreover, we have 
\[\sup_{\{f=r\}} |D^l h_a| \leq O(r^{-\frac{l}{2}})\] 
for each $l \geq 0$. Thus, standard interpolation inequalities imply that 
\[\sup_{\{f=r\}} |Dh_a| \leq O(r^{-\frac{9}{16}}).\] 
On the other hand, we have 
\[\text{\rm div}(h_a) - \frac{1}{2} \, \nabla(\text{\rm tr} \, h_a) = \Delta U_a + \text{\rm Ric}(U_a).\] 
Putting these facts together, obtain 
\[\sup_{\{f=r\}} |\Delta U_a + \text{\rm Ric}(U_a)| \leq O(r^{-\frac{9}{16}}).\] 
Using the estimate $|\text{\rm Ric}(U_a) - D_X U_a| = |[U_a,X]| \leq O(r^{-1})$, we conclude that 
\[\sup_{\{f=r\}} |\Delta U_a + D_X U_a| \leq O(r^{-\frac{9}{16}}).\] 
Finally, the identity 
\[\sum_{a=1}^3 U_a \otimes U_a = r \, (e_1 \otimes e_1 + e_2 \otimes e_2 + O(r^{-\frac{1}{8}}))\] 
follows immediately from (\ref{normalization}). \\

Note that it is enough to define the vector fields $U_1,U_2,U_3$ outside of a compact region. Since we are only interested in the asymptotic behavior near infinity, we can extend the vector fields $U_1,U_2,U_3$ in an arbitrary way into the interior. 

\section{A PDE for the Lie derivative of a vector field}

\label{calculation}

Let us fix a small number $\varepsilon > 0$. For example, $\varepsilon = \frac{1}{100}$ will work. In this section, we consider a vector field $W$ satisfying $\Delta W + D_X W = 0$. Our goal is to derive an elliptic equation for the Lie derivative $\mathscr{L}_W(g)$. 

\begin{theorem} 
\label{pde.for.lie.derivative}
Suppose that $W$ is a vector field satisfying $\Delta W + D_X W = 0$. Then the Lie derivative $\mathscr{L}_W(g)$ satisfies 
\[\Delta_L(\mathscr{L}_W(g)) + \mathscr{L}_X(\mathscr{L}_W(g)) = 0.\] 
\end{theorem} 

\textbf{Proof.} 
Let $g(s)$ be a smooth one-parameter family of metrics with $g(0) = g$. It follows from Proposition 2.3.7 in \cite{Topping} that  
\begin{equation} 
\label{variation.of.ricci}
\frac{\partial}{\partial s} \text{\rm Ric}_{g(s)} \Big |_{s=0} = -\frac{1}{2} \, \Delta_L h + \frac{1}{2} \, \mathscr{L}_Z(g), 
\end{equation}
where $h = \frac{\partial}{\partial s} g(s) \big |_{s=0}$ and 
\[Z = \text{\rm div} \, h - \frac{1}{2} \, \nabla(\text{\rm tr} \, h).\] 
Let us apply the formula (\ref{variation.of.ricci}) to the family of metrics obtained by pulling back $g$ under the one-parameter group of diffeomorphisms generated by $W$. This gives 
\begin{equation} 
\label{lie.derivative.of.ricci}
\mathscr{L}_W(\text{\rm Ric}) = -\frac{1}{2} \, \Delta_L h + \frac{1}{2} \, \mathscr{L}_Z(g), 
\end{equation}
where $h = \mathscr{L}_W(g)$ and 
\[Z = \text{\rm div} \, h - \frac{1}{2} \, \nabla(\text{\rm tr} \, h) = \Delta W + \text{\rm Ric}(W).\] 
Using the relation $\Delta W + D_X W = 0$, we obtain 
\[Z = \Delta W + D_W X = -[X,W].\] 
Substituting this identity into (\ref{lie.derivative.of.ricci}), we conclude that 
\begin{align*} 
\Delta_L(\mathscr{L}_W(g)) 
&= -2 \, \mathscr{L}_W(\text{\rm Ric}) + \mathscr{L}_Z(g) \\ 
&= -\mathscr{L}_W(\mathscr{L}_X(g)) - \mathscr{L}_{[X,W]}(g) \\ 
&= -\mathscr{L}_X(\mathscr{L}_W(g)). 
\end{align*} 
This completes the proof. \\

Applying Theorem \ref{pde.for.lie.derivative} to the vector field $X$ gives the following result:

\begin{proposition}
The vector field $X$ satisfies $\Delta X + D_X X = 0$. Moreover, the Ricci tensor satisfies 
\[\Delta_L(\text{\rm Ric}) + \mathscr{L}_X(\text{\rm Ric}) = 0.\] 
\end{proposition}

\textbf{Proof.} 
Let $h = \mathscr{L}_X(g) = 2 \, \text{\rm Ric}$. The contracted second Bianchi identity implies that 
\[0 = \text{\rm div} \, h - \frac{1}{2} \, \nabla(\text{\rm tr} \, h) = \Delta X + \text{\rm Ric}(X) = \Delta X + D_X X.\] 
Using Theorem \ref{pde.for.lie.derivative}, we obtain 
\[\Delta_L h + \mathscr{L}_X(h) = 0,\] 
as claimed. \\

We note that the identity $\Delta_L(\text{\rm Ric}) + \mathscr{L}_X(\text{\rm Ric}) = 0$ can alternatively be derived from the evolution equation for the Ricci tensor under the Ricci flow (see e.g. \cite{Brendle-book}, Section 2.4).

\section{An elliptic PDE for vector fields}

\label{pde.for.vector.fields}

Throughout this section, we fix a smooth vector field $Q$ on $M$ such that $|Q| \leq O(r^{-\frac{1}{2}-2\varepsilon})$. Our goal is to construct a vector field $V$ on $M$ such that $\Delta V + D_X V = Q$ and $|V| \leq O(r^{\frac{1}{2}-\varepsilon})$. We first establish some auxiliary results.

\begin{lemma}
\label{spherical.harmonics.V}
Let us consider the one-parameter family of shrinking cylinders $(S^2 \times \mathbb{R},\overline{g}(t))$, $t \in (0,1)$, where $\overline{g}(t)$ is given by (\ref{shrinking.cylinders}). Suppose that $\overline{V}(t)$, $t \in (0,1)$, is a one-parameter family of vector fields satisfying the parabolic equation 
\begin{equation} 
\label{evolution.of.vector.field}
\frac{\partial}{\partial t} \overline{V}(t) = \Delta_{\overline{g}(t)} \overline{V}(t) + \text{\rm Ric}_{\overline{g}(t)}(\overline{V}(t)). 
\end{equation} 
Moreover, we assume that $\overline{V}(t)$ is invariant under translations along the axis of the cylinder, and 
\begin{equation} 
\label{boundedness.assumption.V} 
|\overline{V}(t)|_{\overline{g}(t)} \leq 1 
\end{equation} 
for all $t \in (0,\frac{1}{2}]$. Then 
\[\inf_{\lambda \in \mathbb{R}} \sup_{S^2 \times \mathbb{R}} \Big | \overline{V}(t) - \lambda \, \frac{\partial}{\partial z} \Big |_{\overline{g}(t)} \leq L \, (1-t)^{\frac{1}{2}}\] 
for all $t \in [\frac{1}{2},1)$, where $L$ is a positive constant.
\end{lemma} 

\textbf{Proof.} 
Since $\overline{V}(t)$ is invariant under translations along the axis of the cylinder, we may write 
\[\overline{V}(t) = \xi(t) + \eta(t) \, \frac{\partial}{\partial z}\] 
for $t \in (0,1)$, where $\xi(t)$ is a vector field on $S^2$ and $\eta(t)$ is a real-valued function on $S^2$. The parabolic equation (\ref{evolution.of.vector.field}) is equivalent to the following system of equations $\xi(t)$ and $\eta(t)$: 
\begin{align} 
&\frac{\partial}{\partial t} \xi(t) = \frac{1}{2-2t} \, (\Delta_{S^2} \xi(t) + \xi(t)), \label{vec.1} \\ 
&\frac{\partial}{\partial t} \eta(t) = \frac{1}{2-2t} \, \Delta_{S^2} \eta(t). \label{vec.2} 
\end{align}
Moreover, the assumption (\ref{boundedness.assumption.V}) implies 
\begin{align} 
&\sup_{S^2} |\xi(t)|_{g_{S^2}} \leq L_1, \label{est.1} \\ 
&\sup_{S^2} |\eta(t)| \leq L_1 \label{est.2}
\end{align} 
for each $t \in (0,\frac{1}{2}]$, where $L_1$ is a positive constant.

Consider now the operator $\xi \mapsto -\Delta_{S^2} \xi - \xi$, acting on vector fields on $S^2$. It follows from Proposition \ref{laplacian.on.one.forms} that the first eigenvalue of this operator is nonnegative. Using (\ref{vec.1}) and (\ref{est.1}), we conclude that 
\begin{equation} 
\label{est.3} 
\sup_{S^2} |\xi(t)|_{g_{S^2}} \leq L_2
\end{equation} 
for all $t \in [\frac{1}{2},1)$, where $L_2$ is a positive constant. Similarly, using (\ref{vec.2}) and (\ref{est.2}), we can show that 
\begin{equation} 
\label{est.4}
\inf_{\lambda \in \mathbb{R}} \sup_{S^2} |\eta(t) - \lambda| \leq L_3 \, (1-t) 
\end{equation} 
for each $t \in [\frac{1}{2},1)$, where $L_3$ is a positive constant. Combining (\ref{est.3}) and (\ref{est.4}), the assertion follows. \\

\begin{lemma}
\label{max.principle.V}
Let $V$ be a smooth vector field satisfying $\Delta V + D_X V = Q$ in the region $\{f \leq \rho\}$. Then 
\[\sup_{\{f \leq \rho\}} |V| \leq \sup_{\{f=\rho\}} |V| + B \, \rho^{\frac{1}{2}-2\varepsilon}\] 
for some uniform constant $B \geq 1$.
\end{lemma} 

\textbf{Proof.} 
It follows from Kato's inequality that 
\begin{align*} 
\Delta(|V|^2) + \langle X,\nabla(|V|^2) \rangle 
&= 2 \, |DV|^2 + 2 \, \langle V,Q \rangle \\ 
&\geq 2 \, \big | \nabla |V| \big |^2 - 2 \, |Q| \, |V|. 
\end{align*}
This implies 
\[\Delta(|V|) + \langle X,\nabla |V| \rangle \geq -|Q|\] 
when $V \neq 0$. Moreover, using the identity $\Delta f + |\nabla f|^2 = 1$ and the inequality $f \geq 1$, we obtain 
\begin{align*} 
&\Delta(f^{\frac{1}{2}-2\varepsilon}) + \langle X,\nabla(f^{\frac{1}{2}-2\varepsilon}) \rangle \\ 
&= \Big ( \frac{1}{2}-2\varepsilon \Big ) \, f^{-\frac{1}{2}-2\varepsilon} \, (\Delta f + |\nabla f|^2) - \Big ( \frac{1}{4}-4\varepsilon^2 \Big ) \, f^{-\frac{3}{2}-2\varepsilon} \, |\nabla f|^2 \\ 
&\geq \Big ( \frac{1}{2}-2\varepsilon \Big ) \, f^{-\frac{1}{2}-2\varepsilon} - \Big ( \frac{1}{4}-4\varepsilon^2 \Big ) \, f^{-\frac{1}{2}-2\varepsilon} \\ 
&= \Big ( \frac{1}{2}-2\varepsilon \Big )^2 \, f^{-\frac{1}{2}-2\varepsilon}. 
\end{align*}
By assumption, we can find a constant $B \geq 1$ such that 
\[|Q| < \Big ( \frac{1}{2}-2\varepsilon \Big )^2 \, B \, f^{-\frac{1}{2}-2\varepsilon}.\] 
Putting these facts together, we obtain 
\[\Delta (|V| + B \, f^{\frac{1}{2}-2\varepsilon}) + \langle X,\nabla (|V| + B \, f^{\frac{1}{2}-2\varepsilon}) \rangle > 0\] 
when $V \neq 0$. By the maximum principle, the function $|V| + B \, f^{\frac{1}{2}-2\varepsilon}$ attains its maximum on the boundary; that is, 
\[\sup_{\{f \leq \rho\}} (|V| + B \, f^{\frac{1}{2}-2\varepsilon}) \leq \sup_{\{f=\rho\}} |V| + B \, \rho^{\frac{1}{2}-2\varepsilon}.\] From this, the assertion follows. \\

In the following, we consider a sequence of real numbers $\rho_m \to \infty$. Given any integer $m$, there exists a unique vector field $V^{(m)}$ such that 
\[\Delta V^{(m)} + D_X V^{(m)} = Q\] 
in the region $\{f \leq \rho_m\}$ and $V^{(m)} = 0$ on the boundary $\{f=\rho_m\}$. Moreover, we define 
\[A^{(m)}(r) = \inf_{\lambda \in \mathbb{R}} \sup_{\{f=r\}} |V^{(m)} - \lambda \, X|\] 
for $r \leq \rho_m$.

\begin{lemma}
\label{iteration}
Let us fix a real number $\tau \in (0,\frac{1}{2})$ such that $\tau^{-\varepsilon} > 2L$, where $L$ is the constant in Lemma \ref{spherical.harmonics.V}. Then we can find a real number $\rho_0$ and a positive integer $m_0$ such that 
\[2 \, \tau^{-\frac{1}{2}+\varepsilon} \, A^{(m)}(\tau r) \leq A^{(m)}(r) + r^{\frac{1}{2}-\varepsilon}\] 
for all $r \in [\rho_0,\rho_m]$ and all $m \geq m_0$.
\end{lemma}

\textbf{Proof.} 
Suppose that the assertion is false. After passing to a subsequence, we can find a sequence of real numbers $r_m \leq \rho_m$ such that $r_m \to \infty$ and 
\[A^{(m)}(r_m) + r_m^{\frac{1}{2}-\varepsilon} \leq 2 \, \tau^{-\frac{1}{2}+\varepsilon} \, A^{(m)}(\tau r_m)\] 
for all $m$. For each $m$, we choose a real number $\lambda_m$ such that 
\[\sup_{\{f=r_m\}} |V^{(m)} - \lambda_m \, X| = A^{(m)}(r_m).\] 
The vector field $V^{(m)} - \lambda_m \, X$ satisfies the equation  
\[\Delta (V^{(m)} - \lambda_m \, X) + D_X(V^{(m)} - \lambda_m \, X) = Q.\] 
Using Lemma \ref{max.principle.V}, we obtain 
\begin{align*} 
\sup_{\{f \leq r_m\}} |V^{(m)} - \lambda_m \, X| 
&\leq \sup_{\{f=r_m\}} |V^{(m)} - \lambda_m \, X| + B \, r_m^{\frac{1}{2}-2\varepsilon} \\ 
&\leq A^{(m)}(r_m) + r_m^{\frac{1}{2}-\varepsilon} 
\end{align*} 
if $m$ is sufficiently large. Therefore, the vector field 
\[\tilde{V}^{(m)} = \frac{1}{A^{(m)}(r_m) + r_m^{\frac{1}{2}-\varepsilon}} \, (V^{(m)} - \lambda_m \, X)\] 
satisfies 
\begin{equation} 
\label{bound.for.V}
\sup_{\{f \leq r_m\}} |\tilde{V}^{(m)}| \leq 1 
\end{equation} 
if $m$ is sufficiently large. We next define 
\[\hat{g}^{(m)}(t) = r_m^{-1} \, \Phi_{r_m t}^*(g)\] 
and 
\[\hat{V}^{(m)}(t) = r_m^{\frac{1}{2}} \, \Phi_{r_m t}^*(\tilde{V}^{(m)}).\] 
Since $(M,g)$ is a steady Ricci soliton, the metrics $\hat{g}^{(m)}(t)$ form a solution to the Ricci flow. Moreover, the vector fields $\hat{V}^{(m)}(t)$ satisfy the parabolic equation 
\[\frac{\partial}{\partial t} \hat{V}^{(m)}(t) = \Delta_{\hat{g}^{(m)}(t)} \hat{V}^{(m)}(t) + \text{\rm Ric}_{\hat{g}^{(m)}(t)}(\hat{V}^{(m)}(t)) - \hat{Q}^{(m)}(t),\] 
where 
\[\hat{Q}^{(m)}(t) = \frac{r_m^{\frac{3}{2}}}{A^{(m)}(r_m) + r_m^{\frac{1}{2}-\varepsilon}} \, \Phi_{r_mt}^*(Q).\] 
The inequality (\ref{bound.for.V}) implies that 
\[\limsup_{m \to \infty} \sup_{t \in [\delta,1-\delta]} \sup_{\{r_m - \delta^{-1} \, \sqrt{r_m} \leq f \leq r_m + \delta^{-1} \, \sqrt{r_m}\}} |\hat{V}^{(m)}(t)|_{\hat{g}^{(m)}(t)} < \infty\] 
for any given $\delta \in (0,\frac{1}{2})$. Moreover, using the estimate $|Q| \leq O(r^{-\frac{1}{2}-2\varepsilon})$, we obtain 
\[\limsup_{m \to \infty} \sup_{t \in [\delta,1-\delta]} \sup_{\{r_m - \delta^{-1} \, \sqrt{r_m} \leq f \leq r_m + \delta^{-1} \, \sqrt{r_m}\}} |\hat{Q}^{(m)}(t)|_{\hat{g}^{(m)}(t)} = 0\] 
for any given $\delta \in (0,\frac{1}{2})$. 

We now pass to the limit as $m \to \infty$. To that end, we choose a sequence of marked points $p_m \in M$ such that $f(p_m) = r_m$. The sequence $(M,\hat{g}^{(m)}(t),p_m)$ converges in the Cheeger-Gromov sense to a one-parameter family of shrinking cylinders $(S^2 \times \mathbb{R},\overline{g}(t))$, $t \in (0,1)$, where $\overline{g}(t)$ is given by (\ref{shrinking.cylinders}). The rescaled vector fields $r_m^{\frac{1}{2}} \, X$ converge to the axial vector field $\frac{\partial}{\partial z}$ on $S^2 \times \mathbb{R}$. Finally, after passing to a subsequence, the vector fields $\hat{V}^{(m)}(t)$ converge in $C_{loc}^0$ to a one-parameter family of vector fields $\overline{V}(t)$, $t \in (0,1)$, which satisfy the parabolic equation 
\[\frac{\partial}{\partial t} \overline{V}(t) = \Delta_{\overline{g}(t)} \overline{V}(t) + \text{\rm Ric}_{\overline{g}(t)}(\overline{V}(t)).\] 
(The convergence in $C_{loc}^0$ follows from the Arzela-Ascoli theorem together with standard interior estimates for linear parabolic equations; see e.g. \cite{Lieberman}, Theorem 7.22.) Using the identity 
\[\Phi_{\sqrt{r_m} \, s}^*(\hat{V}^{(m)}(t)) = \hat{V}^{(m)} \Big ( t + \frac{s}{\sqrt{r_m}} \Big ),\] 
we conclude that $\Psi_s^*(\overline{V}(t)) = \overline{V}(t)$, where $\Psi_s: S^2 \times \mathbb{R} \to S^2 \times \mathbb{R}$ denotes the flow generated by the axial vector field $-\frac{\partial}{\partial z}$. Hence, $\overline{V}(t)$ is invariant under translations along the axis of the cylinder. Using the estimate (\ref{bound.for.V}), we obtain 
\[|\overline{V}(t)|_{\overline{g}(t)} \leq 1\] 
for all $t \in (0,\frac{1}{2}]$. Using Lemma \ref{spherical.harmonics.V}, we conclude that 
\begin{equation} 
\label{symmetry.V}
\inf_{\lambda \in \mathbb{R}} \sup_{S^2 \times \mathbb{R}} \Big | \overline{V}(t) - \lambda \, \frac{\partial}{\partial z} \Big |_{\overline{g}(t)} \leq L \, (1-t)^{\frac{1}{2}} 
\end{equation}
for all $t \in [\frac{1}{2},1)$. On the other hand, we have 
\begin{align*} 
&\inf_{\lambda \in \mathbb{R}} \sup_{\Phi_{r_m(\tau-1)}(\{f=\tau r_m\})} \Big | \hat{V}^{(m)}(1-\tau) - \lambda \, r_m^{\frac{1}{2}} \, X \Big |_{\hat{g}^{(m)}(1-\tau)} \\ 
&= \inf_{\lambda \in \mathbb{R}} \sup_{\{f=\tau r_m\}} |\tilde{V}^{(m)} - \lambda \, X|_g \\ 
&= \frac{1}{A^{(m)}(r_m)+ r_m^{\frac{1}{2}-\varepsilon}} \, \inf_{\lambda \in \mathbb{R}} \sup_{\{f=\tau r_m\}} |V^{(m)} - \lambda \, X|_g \\ 
&= \frac{A^{(m)}(\tau r_m)}{A^{(m)}(r_m) + r_m^{\frac{1}{2}-\varepsilon}} \\ 
&\geq \frac{1}{2} \, \tau^{\frac{1}{2}-\varepsilon}. 
\end{align*}
Passing to the limit as $m \to \infty$ gives 
\begin{equation} 
\label{asymmetry.V}
\inf_{\lambda \in \mathbb{R}} \sup_{S^2 \times \mathbb{R}} \Big | \overline{V}(1-\tau) - \lambda \, \frac{\partial}{\partial z} \Big |_{\overline{g}(1-\tau)} \geq \frac{1}{2} \, \tau^{\frac{1}{2}-\varepsilon}. 
\end{equation} 
Since $\tau^{-\varepsilon} > 2L$, the inequalities (\ref{symmetry.V}) and (\ref{asymmetry.V}) are in contradiction. This completes the proof of Lemma \ref{iteration}. \\

\begin{proposition} 
\label{uniform.bounds}
There exists a sequence of real numbers $\lambda_m$ such that 
\[\sup_m \sup_{\{f \leq \rho_m\}} f^{-\frac{1}{2}+\varepsilon} \, |V^{(m)} - \lambda_m \, X| < \infty.\] 
\end{proposition}

\textbf{Proof.} 
Let us fix a real number $\tau \in (0,\frac{1}{2})$ so that $\tau^{-\varepsilon} > 2L$, where $L$ is the constant in Lemma \ref{spherical.harmonics.V}. By Lemma \ref{iteration}, we can find a real number $\rho_0$ and a positive integer $m_0$ such that 
\begin{equation} 
\label{iterate}
2 \, \tau^{-\frac{1}{2}+\varepsilon} \, A^{(m)}(\tau r) \leq A^{(m)}(r) + r^{\frac{1}{2}-\varepsilon} 
\end{equation}
for all $r \in [\rho_0,\rho_m]$ and all $m \geq m_0$. Moreover, Lemma \ref{max.principle.V} implies that 
\[\sup_{\rho_0 \leq r \leq \rho_m} A^{(m)}(r) \leq \sup_{\{f \leq \rho_m\}} |V^{(m)}| \leq B \, \rho_m^{\frac{1}{2}-2\varepsilon}.\] 
If we iterate the inequality (\ref{iterate}), we obtain 
\begin{equation} 
\label{key.inequality}
\sup_{m \geq m_0} \sup_{\rho_0 \leq r \leq \rho_m} r^{-\frac{1}{2}+\varepsilon} \, A^{(m)}(r) < \infty. 
\end{equation}
In the next step, we fix a real number $\rho_1 > \rho_0$ such that $\sup_{\{f=\rho_1\}} |X| \geq \frac{1}{2}$. We can find a sequence of real numbers $\lambda_m$ such that 
\[\sup_{\{f=\rho_1\}} |V^{(m)} - \lambda_m \, X| = A^{(m)}(\rho_1)\] 
for each $m$. Applying Lemma \ref{max.principle.V} to the vector field $V^{(m)} - \lambda \, X$, we obtain 
\[\sup_{\{f=\rho_1\}} |V^{(m)} - \lambda \, X| \leq \sup_{\{f=r\}} |V^{(m)} - \lambda \, X| + B \, r^{\frac{1}{2}-2\varepsilon}\] 
for all $r \in [\rho_1,\rho_m]$ and all $\lambda \in \mathbb{R}$. This implies 
\begin{align*} 
&\sup_{\{f=r\}} |V^{(m)} - \lambda_m \, X| \\ 
&\leq \sup_{\{f=r\}} |V^{(m)} - \lambda \, X| + |\lambda - \lambda_m| \\ 
&\leq \sup_{\{f=r\}} |V^{(m)} - \lambda \, X| + 2 \, \sup_{\{f=\rho_1\}} |\lambda \, X - \lambda_m \, X| \\ 
&\leq \sup_{\{f=r\}} |V^{(m)} - \lambda \, X| + 2 \, \sup_{\{f=\rho_1\}} |V^{(m)} - \lambda_m \, X| + 2 \, \sup_{\{f=\rho_1\}} |V^{(m)} - \lambda \, X| \\ 
&\leq 3 \, \sup_{\{f=r\}} |V^{(m)} - \lambda \, X| + 2 \, A^{(m)}(\rho_1) + 2B \, r^{\frac{1}{2}-2\varepsilon} 
\end{align*}  
for all $r \in [\rho_1,\rho_m]$ and all $\lambda \in \mathbb{R}$. Taking the infimum over $\lambda \in \mathbb{R}$ gives 
\[\sup_{\{f=r\}} |V^{(m)} - \lambda_m \, X| \leq 3 \, A^{(m)}(r) + 2 \, A^{(m)}(\rho_1) + 2B \, r^{\frac{1}{2}-2\varepsilon}\] 
for all $r \in [\rho_1,\rho_m]$. Consequently, the inequality (\ref{key.inequality}) implies 
\[\sup_{m \geq m_0} \sup_{\rho_1 \leq r \leq \rho_m} \sup_{\{f=r\}} r^{-\frac{1}{2}+\varepsilon} \, |V^{(m)} - \lambda_m \, X| < \infty,\] 
hence 
\[\sup_{m \geq m_0} \sup_{\{\rho_1 \leq f \leq \rho_m\}} f^{-\frac{1}{2}+\varepsilon} \, |V^{(m)} - \lambda_m \, X| < \infty.\] 
Using Lemma \ref{max.principle.V}, we conclude that 
\[\sup_{m \geq m_0} \sup_{\{f \leq \rho_1\}} |V^{(m)} - \lambda_m \, X| < \infty.\] 
Putting these facts together, the assertion follows. \\

\begin{theorem} 
\label{existence.of.V}
There exists a smooth vector field $V$ such that $\Delta V + D_X V = Q$ and $|V| \leq O(r^{\frac{1}{2}-\varepsilon})$. Moreover, $|DV| \leq O(r^{-\varepsilon})$.
\end{theorem}

\textbf{Proof.}  
By Proposition \ref{uniform.bounds}, we can find a sequence of real numbers $\lambda_m$ such that 
\[\sup_m \sup_{\{f \leq \rho_m\}} f^{-\frac{1}{2}+\varepsilon} \, |V^{(m)} - \lambda_m \, X| < \infty.\] 
Moreover, the vector field $V^{(m)} - \lambda_m \, X$ solves the equation 
\[\Delta (V^{(m)} - \lambda_m \, X) + D_X(V^{(m)} - \lambda_m \, X) = Q\] 
in the region $\{f \leq \rho_m\}$. Hence, after passing to a subsequence if necessary, the vector fields $V^{(m)} - \lambda_m \, X$ converge to a smooth vector field $V$ satisfying $\Delta V + D_X V = Q$ and $|V| \leq O(r^{\frac{1}{2}-\varepsilon})$.

It remains to show that $|DV| \leq O(r^{-\varepsilon})$. In order to prove this, we use the standard interior regularity theory for parabolic equations. Consider a sequence $r_m \to \infty$, and let 
\[\hat{g}^{(m)}(t) = r_m^{-1} \, \Phi_{r_m t}^*(g)\] 
for $t \in [-\frac{1}{2},0]$. Moreover, we define 
\[\hat{V}^{(m)}(t) = \Phi_{r_m t}^*(V)\] 
and 
\[\hat{Q}^{(m)}(t) = r_m \, \Phi_{r_mt}^*(Q)\] 
for $t \in [-\frac{1}{2},0]$. The vector fields $\hat{V}^{(m)}(t)$ satisfy the parabolic equation 
\[\frac{\partial}{\partial t} \hat{V}^{(m)}(t) = \Delta_{\hat{g}^{(m)}(t)} \hat{V}^{(m)}(t) + \text{\rm Ric}_{\hat{g}^{(m)}(t)}(\hat{V}^{(m)}(t)) - \hat{Q}^{(m)}(t).\] 
Moreover, since $|Q| \leq O(r^{-\frac{1}{2}-2\varepsilon})$, we have  
\[\sup_{t \in [-\frac{1}{2},0]} \sup_{\{r_m - \sqrt{r_m} \leq f \leq r_m + \sqrt{r_m}\}} |\hat{Q}^{(m)}(t)|_{\hat{g}^{(m)}(t)} \leq O(r_m^{-2\varepsilon}).\] 
Using standard interior estimates for parabolic equations, we obtain 
\begin{align*} 
\sup_{\{f=r_m\}} |D\hat{V}^{(m)}(0)|_{\hat{g}^{(m)}(0)} 
&\leq C \, \sup_{t \in [-\frac{1}{2},0]} \sup_{\{r_m - \sqrt{r_m} \leq f \leq r_m + \sqrt{r_m}\}} |\hat{V}^{(m)}(t)|_{\hat{g}^{(m)}(t)} \\ 
&+ C \, \sup_{t \in [-\frac{1}{2},0]} \sup_{\{r_m - \sqrt{r_m} \leq f \leq r_m + \sqrt{r_m}\}} |\hat{Q}^{(m)}(t)|_{\hat{g}^{(m)}(t)} \\ 
&\leq O(r_m^{-\varepsilon}). 
\end{align*} From this, we deduce that 
\[\sup_{\{f=r_m\}} |DV| \leq O(r_m^{-\varepsilon}),\] 
as claimed.

\section{Analysis of the Lichnerowicz equation}

\label{analysis.of.lich.eq}

\begin{lemma}
\label{spherical.harmonics.h}
Let us consider the shrinking cylinders $(S^2 \times \mathbb{R},\overline{g}(t))$, $t \in (0,1)$, where $\overline{g}(t)$ is given by (\ref{shrinking.cylinders}). Suppose that $\overline{h}(t)$, $t \in (0,1)$, is a one-parameter family of $(0,2)$-tensors satisfying the parabolic Lichnerowicz equation 
\begin{equation} 
\label{parabolic.lichnerowicz.equation}
\frac{\partial}{\partial t} \overline{h}(t) = \Delta_{L,\overline{g}(t)} \overline{h}(t). 
\end{equation} 
Moreover, we assume that $\overline{h}(t)$ is invariant under translations along the axis of the cylinder, and 
\begin{equation} 
\label{boundedness.assumption.h}
|\overline{h}(t)|_{\overline{g}(t)} \leq (1-t)^{-1} 
\end{equation}
for all $t \in (0,\frac{1}{2}]$. Then we have 
\[\inf_{\lambda \in \mathbb{R}} \sup_{S^2 \times \mathbb{R}} \big | \overline{h}(t) - \lambda \, \text{\rm Ric}_{\overline{g}(t)} \big |_{\overline{g}(t)} \leq N\] 
for all $t \in [\frac{1}{2},1)$, where $N$ is a positive constant.
\end{lemma}

\textbf{Proof.} 
Since $\overline{h}(t)$ is invariant under translations along the axis of the cylinder, we may write 
\[\overline{h}(t) = \chi(t) + dz \otimes \sigma(t) + \sigma(t) \otimes dz + \beta(t) \, dz \otimes dz\] 
for $t \in (0,1)$, where $\chi(t)$ is a symmetric $(0,2)$ tensor on $S^2$, $\sigma(t)$ is a one-form on $S^2$, and $\beta(t)$ is a real-valued function on $S^2$. The parabolic Lichnerowicz equation (\ref{parabolic.lichnerowicz.equation}) is equivalent to the following system of equations for $\chi(t)$, $\sigma(t)$, and $\beta(t)$: 
\begin{align} 
&\frac{\partial}{\partial t} \chi(t) = \frac{1}{2-2t} \, (\Delta_{S^2} \chi(t) - 4 \, \overset{\text{\rm o}}{\chi}(t)), \label{lich.1} \\ 
&\frac{\partial}{\partial t} \sigma(t) = \frac{1}{2-2t} \, (\Delta_{S^2} \sigma(t) - \sigma(t)), \label{lich.2} \\ 
&\frac{\partial}{\partial t} \beta(t) = \frac{1}{2-2t} \, \Delta_{S^2} \beta(t). \label{lich.3}
\end{align}
Here, $\overset{\text{\rm o}}{\chi}(t)$ denotes the trace-free part of $\chi(t)$ with respect to the standard metric on $S^2$. Moreover, the assumption (\ref{boundedness.assumption.h}) implies 
\begin{align} 
&\sup_{S^2} |\chi(t)|_{g_{S^2}} \leq N_1, \label{estimate.1} \\ 
&\sup_{S^2} |\sigma(t)|_{g_{S^2}} \leq N_1, \label{estimate.2} \\ 
&\sup_{S^2} |\beta(t)| \leq N_1 \label{estimate.3}
\end{align} 
for each $t \in (0,\frac{1}{2}]$, where $N_1$ is a positive constant.

Let us consider the operator $\chi \mapsto -\Delta_{S^2} \chi + 4 \, \overset{\text{\rm o}}{\chi}$, acting on symmetric $(0,2)$-tensors on $S^2$. The first eigenvalue of this operator is equal to $0$, and the associated eigenspace is spanned by $g_{S^2}$. Moreover, all other eigenvalues are at least $2$ (cf. Proposition \ref{laplacian.on.tensors} below). Hence, it follows from (\ref{lich.1}) and (\ref{estimate.1}) that 
\begin{equation} 
\label{estimate.4}
\inf_{\lambda \in \mathbb{R}} \sup_{S^2} |\chi(t) - \lambda \, g_{S^2}|_{g_{S^2}} \leq N_2 \, (1-t) 
\end{equation} 
for all $t \in [\frac{1}{2},1)$, where $N_2$ is a positive constant. We next consider the operator $\sigma \mapsto -\Delta_{S^2} \sigma + \sigma$, acting on one-forms on $S^2$. By Proposition \ref{laplacian.on.one.forms}, the first eigenvalue of this operator is at least $2$. Using (\ref{lich.2}) and (\ref{estimate.2}), we deduce that 
\begin{equation}
\label{estimate.5}
\sup_{S^2} |\sigma(t)|_{g_{S^2}} \leq N_3 \, (1-t) 
\end{equation}
for all $t \in [\frac{1}{2},1)$, where $N_3$ is a positive constant. Finally, using (\ref{lich.3}) and (\ref{estimate.3}), we obtain 
\begin{equation} 
\label{estimate.6}
\sup_{S^2} |\beta(t)| \leq N_4 
\end{equation} 
for all $t \in [\frac{1}{2},1)$, where $N_4$ is a positive constant. Combining (\ref{estimate.4}), (\ref{estimate.5}), and (\ref{estimate.6}), the assertion follows. \\

In the following, we study the equation $\Delta_L h + \mathscr{L}_X(h) = 0$ on $(M,g)$.

\begin{lemma}
\label{max.principle.h}
Let $h$ be a solution of the Lichnerowicz-type equation 
\[\Delta_L h + \mathscr{L}_X(h) = 0\] 
on the region $\{f \leq \rho\}$. Then 
\[\sup_{\{f \leq \rho\}} f \, |h| \leq B \, \rho \sup_{\{f=\rho\}} |h|,\] 
where $B$ is a positive constant that does not depend on $\rho$.
\end{lemma}

\textbf{Proof.} 
By a result of Anderson and Chow \cite{Anderson-Chow}, we have 
\[\Delta \Big ( \frac{|h|^2}{R^2} \Big ) + \Big \langle X + 2 \, \frac{\nabla R}{R},\nabla \Big ( \frac{|h|^2}{R^2} \Big ) \Big \rangle \geq 0.\] 
Applying the maximum principle, we obtain 
\[\sup_{\{f \leq \rho\}} \frac{|h|}{R} \leq \sup_{\{f=\rho\}} \frac{|h|}{R}.\] 
Since $\sup_M f \, R < \infty$ and $\inf_M f \, R > 0$, the assertion follows. \\

\begin{theorem}
\label{lichnerowicz.equation}
Let $h$ be a solution of the Lichnerowicz-type equation 
\[\Delta_L h + \mathscr{L}_X(h) = 0\] 
such that $|h| \leq O(r^{-\varepsilon})$. Then $h = \lambda \, \text{\rm Ric}$ for some constant $\lambda \in \mathbb{R}$.
\end{theorem} 

\textbf{Proof.} Let 
\[A(r) = \inf_{\lambda \in \mathbb{R}} \sup_{\{f=r\}} |h - \lambda \, \text{\rm Ric}|.\] 
Clearly, $A(r) \leq \sup_{\{f=r\}} |h| \leq O(r^{-\varepsilon})$. We consider two cases: 

\textit{Case 1:} Suppose that there exists a sequence of real numbers $r_m \to \infty$ such that $A(r_m) = 0$ for all $m$. For each $m$, we choose a real number $\lambda_m$ such that 
\[\sup_{\{f=r_m\}} |h - \lambda_m \, \text{\rm Ric}| = A(r_m) = 0.\] 
Applying Lemma \ref{max.principle.h} to the tensor $h - \lambda_m \, \text{\rm Ric}$, we obtain 
\[\sup_{\{f \leq r_m\}} f \, |h - \lambda_m \, \text{\rm Ric}| \leq B \, r_m \sup_{\{f=r_m\}} |h - \lambda_m \, \text{\rm Ric}| = 0.\] 
Therefore, we have $h - \lambda_m \, \text{\rm Ric} = 0$ in the region $\{f \leq r_m\}$. Consequently, the sequence $\lambda_m$ is constant and $h$ is a constant multiple of the Ricci tensor.

\textit{Case 2:} Suppose now that $A(r) > 0$ when $r$ is sufficiently large. We fix a real number $\tau \in (0,\frac{1}{2})$ such that $\tau^{-\varepsilon} > 2N \, B$, where $N$ is the constant in Lemma \ref{spherical.harmonics.h} and $B$ is the constant in Lemma \ref{max.principle.h}. Since $A(r) \leq O(r^{-\varepsilon})$, we can find a sequence of real numbers $r_m \to \infty$ such that 
\[A(r_m) \leq 2 \, \tau^\varepsilon \, A(\tau r_m)\] 
for all $m$. For each $m$, we choose a real number $\lambda_m$ such that 
\[\sup_{\{f=r_m\}} |h - \lambda_m \, \text{\rm Ric}| = A(r_m).\] 
The tensor 
\[\tilde{h}^{(m)} = \frac{1}{A(r_m)} \, (h - \lambda_m \, \text{\rm Ric})\] 
satisfies the Lichnerowicz-type equation 
\[\Delta_L \tilde{h}^{(m)} + \mathscr{L}_X(\tilde{h}^{(m)}) = 0.\] 
Using Lemma \ref{max.principle.h}, we obtain 
\begin{equation} 
\label{bound.for.h}
\sup_{\{f=r\}} |\tilde{h}^{(m)}| \leq \frac{B \, r_m}{r} \, \sup_{\{f=r_m\}} |\tilde{h}^{(m)}| = \frac{B \, r_m}{r \, A(r_m)} \, \sup_{\{f=r_m\}} |h - \lambda_m \, \text{\rm Ric}| = \frac{B \, r_m}{r} 
\end{equation}
for $r \leq r_m$. 

We now define 
\[\hat{g}^{(m)}(t) = r_m^{-1} \, \Phi_{r_m t}^*(g)\] 
and 
\[\hat{h}^{(m)}(t) = r_m^{-1} \, \Phi_{r_m t}^*(\tilde{h}^{(m)}).\] 
Since $(M,g)$ is a steady Ricci soliton, the metrics $\hat{g}^{(m)}(t)$ evolve by the Ricci flow. Moreover, the tensors $\hat{h}^{(m)}(t)$ satisfy the parabolic Lichnerowicz equation 
\[\frac{\partial}{\partial t} \hat{h}^{(m)}(t) = \Delta_{L,\hat{g}^{(m)}(t)} \hat{h}^{(m)}(t).\] 
It follows from (\ref{bound.for.h}) that 
\[\limsup_{m \to \infty} \sup_{t \in [\delta,1-\delta]} \sup_{\{r_m - \delta^{-1} \, \sqrt{r_m} \leq f \leq r_m + \delta^{-1} \, \sqrt{r_m}\}} |\hat{h}^{(m)}(t)|_{\hat{g}^{(m)}(t)} < \infty\] 
for any given $\delta \in (0,\frac{1}{2})$. 

We next take the limit as $m \to \infty$. As above, we choose a sequence of marked points $p_m \in M$ such that $f(p_m) = r_m$. The sequence $(M,\hat{g}^{(m)}(t),p_m)$ converges in the Cheeger-Gromov sense to a one-parameter family of shrinking cylinders $(S^2 \times \mathbb{R},\overline{g}(t))$, $t \in (0,1)$, where $\overline{g}(t)$ is given by (\ref{shrinking.cylinders}). The rescaled vector fields $r_m^{\frac{1}{2}} \, X$ converge to the axial vector field $\frac{\partial}{\partial z}$ on $S^2 \times \mathbb{R}$. Finally, after passing to a subsequence, the tensors $\hat{h}^{(m)}(t)$ converge in $C_{loc}^\infty$ to a one-parameter family of tensor fields $\overline{h}(t)$, $t \in (0,1)$, which satisfy the parabolic Lichnerowicz equation 
\[\frac{\partial}{\partial t} \overline{h}(t) = \Delta_{L,\overline{g}(t)} \overline{h}(t).\] 
Using the identity 
\[\Phi_{\sqrt{r_m} \, s}^*(\hat{h}^{(m)}(t)) = \hat{h}^{(m)} \Big ( t + \frac{s}{\sqrt{r_m}} \Big ),\] 
we obtain $\Psi_s^*(\overline{h}(t)) = \overline{h}(t)$, where $\Psi_s: S^2 \times \mathbb{R} \to S^2 \times \mathbb{R}$ denotes the flow generated by the axial vector field $-\frac{\partial}{\partial z}$. In other words, $\overline{h}(t)$ is invariant under translations along the axis of the cylinder. Moreover, the estimate (\ref{bound.for.h}) implies 
\[|\overline{h}(t)|_{\overline{g}(t)} \leq B \, (1-t)^{-1}\] 
for all $t \in (0,\frac{1}{2}]$. Using Lemma \ref{spherical.harmonics.h}, we conclude that 
\begin{equation} 
\label{symmetry.h}
\inf_{\lambda \in \mathbb{R}} \sup_{S^2 \times \mathbb{R}} \big | \overline{h}(t) - \lambda \, \text{\rm Ric}_{\overline{g}(t)} \big |_{\overline{g}(t)} \leq N \, B 
\end{equation}
for all $t \in [\frac{1}{2},1)$. On the other hand, we have 
\begin{align*} 
&\inf_{\lambda \in \mathbb{R}} \sup_{\Phi_{r_m(\tau-1)}(\{f=\tau r_m\})} \Big | \hat{h}^{(m)}(1-\tau) - \lambda \, \text{\rm Ric}_{\hat{g}^{(m)}(1-\tau)} \Big |_{\hat{g}^{(m)}(1-\tau)} \\ 
&= \inf_{\lambda \in \mathbb{R}} \sup_{\{f=\tau r_m\}} |\tilde{h}^{(m)} - \lambda \, \text{\rm Ric}_g|_g \\ 
&= \frac{1}{A(r_m)} \, \inf_{\lambda \in \mathbb{R}} \sup_{\{f=\tau r_m\}} |h - \lambda \, \text{\rm Ric}_g|_g \\ 
&= \frac{A(\tau r_m)}{A(r_m)} \\ 
&\geq \frac{1}{2} \, \tau^{-\varepsilon}. 
\end{align*}
Taking the limit as $m \to \infty$ gives 
\begin{equation} 
\label{asymmetry.h}
\inf_{\lambda \in \mathbb{R}} \sup_{S^2 \times \mathbb{R}} \big | \overline{h}(1-\tau) - \lambda \, \text{\rm Ric}_{\overline{g}(1-\tau)} \big |_{\overline{g}(1-\tau)} \geq \frac{1}{2} \, \tau^{-\varepsilon}. 
\end{equation}
Since $\tau^{-\varepsilon} > 2N \, B$, the inequalities (\ref{symmetry.h}) and (\ref{asymmetry.h}) are in contradiction. This completes the proof of Theorem \ref{lichnerowicz.equation}. \\

\section{Proof of Theorem \ref{main.theorem}}

\label{sym.prin}

Combining Theorems \ref{pde.for.lie.derivative}, \ref{existence.of.V}, and \ref{lichnerowicz.equation}, we obtain the following symmetry principle: 

\begin{theorem}
\label{symmetry.principle}
Suppose that $U$ is a vector field on $(M,g)$ such that $|\mathscr{L}_U(g)| \leq O(r^{-2\varepsilon})$ and $|\Delta U + D_X U| \leq O(r^{-\frac{1}{2}-2\varepsilon})$ for some small constant $\varepsilon>0$. Then there exists a vector field $\hat{U}$ on $(M,g)$ such that $\mathscr{L}_{\hat{U}}(g) = 0$, $[\hat{U},X] = 0$, $\langle \hat{U},X \rangle = 0$, and $|\hat{U} - U| \leq O(r^{\frac{1}{2}-\varepsilon})$.
\end{theorem}

\textbf{Proof.} 
By Theorem \ref{existence.of.V}, we can find a smooth vector field $V$ such that 
\[\Delta V + D_X V = \Delta U + D_X U\] 
and $|V| \leq O(r^{\frac{1}{2}-\varepsilon})$. Moreover, the covariant derivative of $V$ satisfies $|DV| \leq O(r^{-\varepsilon})$. We now define $W = U-V$ and $h = \mathscr{L}_W(g)$. Since $W$ satisfies the equation $\Delta W + D_X W = 0$, Theorem \ref{pde.for.lie.derivative} implies that the tensor $h$ satisfies the Lichnerowicz-type equation 
\[\Delta_L h + \mathscr{L}_X(h) = 0.\] 
Moreover, $|h| \leq O(r^{-\varepsilon})$. Hence, it follows from Theorem \ref{lichnerowicz.equation} that $h = \lambda \, \text{\rm Ric}$ for some constant $\lambda \in \mathbb{R}$. Therefore, the vector field $\hat{U} := U - V - \frac{1}{2} \, \lambda \, X$ is a Killing vector field. The relation $\mathscr{L}_{\hat{U}}(g) = 0$ implies that $\Delta \hat{U} + \text{\rm Ric}(\hat{U}) = 0$. On the other hand, we have $\Delta \hat{U} + D_X \hat{U} = 0$ by definition of $V$. Thus, we conclude that $[\hat{U},X] = \text{\rm Ric}(\hat{U}) - D_X \hat{U} = 0$. Finally, since $\hat{U}$ is a Killing vector field, we have 
\[D^2(\mathscr{L}_{\hat{U}}(f)) = \mathscr{L}_{\hat{U}}(D^2 f) = \frac{1}{2} \, \mathscr{L}_{\hat{U}}(\mathscr{L}_X(g)) = \frac{1}{2} \, \mathscr{L}_X(\mathscr{L}_{\hat{U}}(g)) = 0.\] 
Consequently, the function $\mathscr{L}_{\hat{U}}(f) = \langle \hat{U},X \rangle$ is constant. Since $X$ vanishes at the point where $f$ attains its minimum, we conclude that the function $\langle \hat{U},X \rangle$ vanishes identically. This completes the proof of Theorem \ref{symmetry.principle}. \\

If we apply Theorem \ref{symmetry.principle} to the vector fields $U_1,U_2,U_3$ constructed in Proposition \ref{rotation.vector.fields}, we can draw the following conclusion: 

\begin{corollary}
\label{rotation.vector.fields.2}
We can find vector fields $\hat{U}_1,\hat{U}_2,\hat{U}_3$ on $(M,g)$ such that $\mathscr{L}_{\hat{U}_a}(g) = 0$, $[\hat{U}_a,X] = 0$, and $\langle \hat{U}_a,X \rangle = 0$. Moreover, we have 
\[\sum_{a=1}^3 \hat{U}_a \otimes \hat{U}_a = r \, (e_1 \otimes e_1 + e_2 \otimes e_2 + O(r^{-\varepsilon})),\] 
where $\{e_1,e_2\}$ is a local orthonormal frame on the level set $\{f=r\}$. 
\end{corollary}

In particular, we have $\text{\rm span}\{\hat{U}_1,\hat{U}_2,\hat{U}_3\} = \text{\rm span}\{e_1,e_2\}$ at each point in $M$ which is sufficiently far out near infinity. This shows that $(M,g)$ is exactly rotationally symmetric near infinity. From this, Theorem \ref{main.theorem} follows easily.

\appendix 

\section{The eigenvalues of some elliptic operators on $S^2$}

In this section, we collect some well-known results concerning the eigenvalues of certain elliptic operators on $S^2$. In the following, $g_{S^2}$ will denote the standard metric on $S^2$ with constant Gaussian curvature $1$. 

\begin{proposition} 
\label{laplacian.on.one.forms}
Let $\sigma$ be a one-form on $S^2$ satisfying 
\[\Delta_{S^2} \sigma + \mu \, \sigma = 0,\] 
where $\Delta_{S^2}$ denotes the rough Laplacian and $\mu \in (-\infty,1)$ is a constant. Then $\sigma = 0$. 
\end{proposition} 

\textbf{Proof.} 
We can find a real-valued function $\alpha$ and a two-form $\omega$ such that $\sigma = d\alpha + d^*\omega$. Using the Bochner formula for one-forms, we obtain 
\begin{align*} 
0 
&= \Delta_{S^2} \sigma + \mu \, \sigma \\ 
&= -dd^*\sigma - d^*d\sigma + (\mu+1) \, \sigma \\ 
&= -dd^*d\alpha - d^*dd^*\omega + (\mu+1) \, (d\alpha + d^* \omega) \\ 
&= d(\Delta_{S^2} \alpha + (\mu+1) \, \alpha) + d^*(\Delta_{S^2} \omega + (\mu+1) \, \omega). 
\end{align*} 
Consequently, the function $\Delta_{S^2} \alpha + (\mu+1) \, \alpha$ is constant, and the two-form $\Delta_{S^2} \omega + (\mu+1) \, \omega$ is a constant multiple of the volume form. Since $\mu+1 < 2$, we conclude that $\alpha$ is constant and $\omega$ is a constant multiple of the volume form. Thus, $\sigma = 0$, as claimed. \\

\begin{proposition} 
\label{laplacian.on.tensors}
Let $\chi$ be a symmetric $(0,2)$-tensor on $S^2$ satisfying 
\[\Delta_{S^2} \chi - 4 \, \overset{\text{\rm o}}{\chi} + \mu \, \chi = 0,\] 
where $\overset{\text{\rm o}}{\chi}$ denotes the trace-free part of $\chi$ and $\mu \in (-\infty,2)$ is a constant. Then $\chi$ is a constant multiple of $g_{S^2}$.
\end{proposition} 

\textbf{Proof.} 
The trace of $\chi$ satisfies 
\[\Delta_{S^2}(\text{\rm tr} \, \chi) + \mu \, (\text{\rm tr} \, \chi) = 0.\] 
Since $\mu < 2$, we conclude that $\text{\rm tr} \, \chi$ is constant. Moreover, the trace-free part of $\chi$ satisfies 
\[\Delta_{S^2} \overset{\text{\rm o}}{\chi} + (\mu - 4) \, \overset{\text{\rm o}}{\chi} = 0.\] 
Since $\mu-4 < 0$, it follows that $\overset{\text{\rm o}}{\chi} = 0$. Putting these facts together, the assertion follows.

\end{document}